\documentclass{amsart}

\usepackage{amssymb}
\usepackage{amscd}
\usepackage{epsfig}

\address{Graduate School of Mathematics, Nagoya University, Furo-cho, Chikusa-ku, Nagoya, 464-8602, Japan}

\email{furusho@math.nagoya-u.ac.jp}

\thanks{{\bf Version of August. 16th.2011}}
\newtheorem{thm}{Theorem}%[section]

\newtheorem{prop}[thm]{Proposition}  

\theoremstyle{remark}

\theoremstyle{definition}
\newtheorem{thm-defn}[thm]{Theorem-Definition} 
\newtheorem{defn}[thm]{Definition}
\newtheorem{rem}[thm]{Remark}

\newtheorem{eg}[thm]{Examples}       

\newtheorem{q}[thm]{Question}

\begin{document}
\bibliographystyle{amsalpha+}
%\maketitle
%%%%%%%%%%%%%%%%%%%%%%%%%%%%%%%%%%%%%%%%%%%%%%%%%%%%%%%%%%%%%%%%%%%%%%
\begin{center}
{\Large{\bf{Four Groups Related to Associators}}}
\end{center}

\bigskip

\begin{center}
{\large Hidekazu Furusho}\\
%Graduate School of Mathematics,Nagoya University
\end{center}
\bigskip
\begin{center}
{\tt Mathematische Arbeitstagung 24th June-1st July.2011.}
\end{center}
\bigskip

\begin{flushright}
{\it Dedicated to Professor Don Zagier \\
on the occasion of his 60th birthday}
\end{flushright}

\bigskip
\noindent

%%%%%%%%%%%%%%%%%%%%%%%%%%%%%%%%%%%%%%%%%%%%%%%%%%%%%%%%%%%%%%%%%%%%%%
\section{Associators}
We recall the definition of associators \cite{Dr} and explain our main results
in \cite{F10,F11}
which are on the defining equations of associators.

Let us fix notations:
Let $k$ be a field of characteristic $0$ and
$\bar k$ its algebraic closure.
Denote by
$U\frak F_2=k\langle\langle X_0,X_1\rangle\rangle$
a non-commutative formal power series ring,
%which is 
a universal enveloping algebra of the completed free Lie algebra
$\frak F_2$ with two variables $X_0$ and $X_1$.
Its element $\varphi=\varphi(X_0,X_1)$ is called {\it group-like}
\footnote{
It is equivalent to $\varphi\in\exp{\frak F}_2$.
}
if it satisfies
\begin{equation}\label{shuffle}
\Delta (\varphi)=\varphi\otimes \varphi
\text{ and } 
\varphi(0,0)=1
\end{equation}
with $\Delta(X_0)=X_0\otimes 1+1\otimes X_0$ and
$\Delta(X_1)=X_1\otimes 1+1\otimes X_1$.
% and its constant term is equal to $1$.
For any $k$-algebra homomorphism $\iota:U\frak F_2\to S$,
the image $\iota(\varphi)\in S$ is denoted 
by $\varphi(\iota(X_0),\iota(X_1))$.

\begin{defn}[\cite{Dr}]
A pair $(\mu,\varphi)$ with a {\it non-zero} element $\mu$ in $k$ and a group-like series $\varphi=\varphi(X_0,X_1)\in U\frak F_2$ is called an {\it associator}
if it satisfies {\it one pentagon equation} in $U\frak a_4$
\begin{equation}\label{pentagon}
\varphi(t_{12},t_{23}+t_{24})
\varphi(t_{13}+t_{23},t_{34})=
\varphi(t_{23},t_{34})
\varphi(t_{12}+t_{13},t_{24}+t_{34})
\varphi(t_{12},t_{23})
\end{equation}
and
{\it two hexagon equations} in $U\frak a_3$
\begin{equation}\label{hexagon}
\exp\{\frac{\mu (t_{13}+t_{23})}{2}\}=
\varphi(t_{13},t_{12})\exp\{\frac{\mu t_{13}}{2}\}
\varphi(t_{13},t_{23})^{-1}
\exp\{\frac{\mu t_{23}}{2}\} \varphi(t_{12},t_{23}),
\end{equation}
\begin{equation}\label{hexagon-b}
\exp\{\frac{\mu (t_{12}+t_{13})}{2}\}=
\varphi(t_{23},t_{13})^{-1}\exp\{\frac{\mu t_{13}}{2}\}
\varphi(t_{12},t_{13})
\exp\{\frac{\mu t_{12}}{2}\} \varphi(t_{12},t_{23})^{-1}.
\end{equation}
\end{defn}

Here $U\frak a_3$ (resp. $U\frak a_4$) means the universal enveloping algebra of
the {\it completed pure braid Lie algebra} $\frak a_3$ (resp. $\frak a_4$) 
over $k$  with 3 (resp. 4) strings,
generated by $t_{ij}$ ($1\leqslant i,j \leqslant 3$ (resp. $4$)) 
with defining relations 
$$t_{ii}=0, \  t_{ij}=t_{ji}, \ [t_{ij},t_{ik}+t_{jk}]=0
\text{ ($i$,$j$,$k$: all distinct)}$$
$$\text{and  }\ [t_{ij},t_{kl}]=0 
\text{ ($i$,$j$,$k$,$l$: all distinct).}$$

\begin{rem}
(i). 
Drinfeld \cite{Dr} proved that such a pair always exists for any filed $k$ of
characteristic $0$.

(ii).
The equations \eqref{pentagon}$\sim$\eqref{hexagon-b} reflect the three axioms of 
braided monoidal categories %introduced by Joyal and Street 
\cite{JS}.
We note that for any $k$-linear {\it infinitesimal} tensor category $\mathcal C$ 
each associator  gives a structure of 
braided  %(actually ribbon) 
monoidal category on 
${\mathcal C}[[h]]$ (cf.\cite{C, Dr}).
Here ${\mathcal C}[[h]]$ means the category whose set of objects is equal to that of
${\mathcal C}$ and whose set of morphism $Mor_{{\mathcal C}[[h]]}(X,Y)$ is
$Mor_{\mathcal C}(X,Y)\otimes k[[h]]$ ($h$: a parameter).

(iii).
Associators are essential for construction of quasi-triangular quasi-Hopf
quantized universal enveloping algebras \cite{Dr}.

(iv).
Le and Murakami \cite{LM} and Bar-Natan \cite{BN} gave a reconstruction of
universal Vassiliev knot invariant (Kontsevich invariant \cite{K,Ba95})
in a combinatorial way by using associators.
\end{rem}

Our first result is the implication of two hexagon equations from
one pentagon equation.

\begin{thm}[\cite{F10}]\label{pentagon-hexagons}
Let $\varphi=\varphi(X_0,X_1)$ be a %commutator 
group-like element of $U\frak F_2$.
Suppose that $\varphi$ satisfies the pentagon equation \eqref{pentagon}.
Then there exists $\mu\in\bar k$ (unique up to signature)
such that the pair $(\mu,\varphi)$ satisfies
two hexagon equations \eqref{hexagon} and \eqref{hexagon-b}.
\end{thm}

Recently several different proofs of the above theorem were obtained
(see \cite{AT, BD, W}).

One of the nice examples of associators is the Drinfeld associator below.

\begin{eg}\label{Drinfeld associator}
The {\it Drinfeld associator} $\varPhi_{KZ}=\varPhi_{KZ}(X_0,X_1)\in
{\bold C}\langle\langle X_0,X_1\rangle\rangle$
is
defined to be the quotient $\varPhi_{KZ}=G_1(z)^{-1}G_0(z)$
where $G_0$ and $G_1$ are the solutions
of the {\it formal KZ (Knizhnik-Zamolodchikov) equation},
the following differential equation over ${\bf C}\backslash\{0,1\}$
with $G(z)$ valued on ${\bold C}\langle\langle X_0,X_1\rangle\rangle$
$$
\frac{d}{dz}G(z)=\bigl(\frac{X_0}{z}+\frac{X_1}{z-1}\bigr)G(z),
$$
such that $G_0(z)\approx z^{X_0}$ when $z\to 0$ and
$G_1(z)\approx (1-z)^{X_1}$ when $z\to 1$  (cf.\cite{Dr}).
It is shown in \cite{Dr} that the pair
$(2\pi\sqrt{-1},\varPhi_{KZ})$
forms an associator for $k=\bf C$.
Namely $\varPhi_{KZ}$
satisfies \eqref{shuffle}$\sim$\eqref{hexagon-b} with $\mu=2\pi\sqrt{-1}$.
\end{eg}

\begin{rem}
(i).
The Drinfeld associator is expressed as follows:
$$
\varPhi_{KZ}(X_0,X_1)=1+
\sum
(-1)^m\zeta(k_1,\cdots,k_m)X_0^{k_m-1}X_1\cdots X_0^{k_1-1}X_1+
\text{(regularized terms)}.
$$
Here $\zeta(k_1,\cdots,k_m)$ is the
{\it multiple zeta value} (MZV in short),
the real number defined by the following power series
\begin{equation}\label{MZV}
\zeta(k_1,\cdots,k_m)
:=\sum_{0<n_1<\cdots<n_m}\frac{1}
%{\zeta_1^{n_1}\cdots \zeta_m^{n_m}}
{n_1^{k_1}\cdots n_m^{k_m}}
\end{equation}
for $m$, $k_1$,\dots, $k_m\in {\bf N} (={\bf Z}_{>0})$
with $k_m>1$ (its convergent condition).
All the coefficients of $\varPhi_{KZ}$ including its regularized terms
are explicitly calculated in terms of
MZV's  in \cite{F03} proposition 3.2.3
by Le-Murakami's method in \cite{LMb}.

(ii).
Since MZV's are coefficients of $\varPhi_{KZ}$,
the equations \eqref{shuffle}$\sim$\eqref{hexagon-b} for
$(\mu,\varphi)=(2\pi\sqrt{-1},\varPhi_{KZ})$ yield
algebraic relations among them,
which are called {\it associator relations}.
It is expected that the associator relations might
produce all algebraic relations among MZV's.

\end{rem}

Various relations among MZV's have been found and studied so far.
The regularised double shuffle relations which were
initially introduced by Zagier and Ecalle in early 90's
might be one of the most fascinating ones. 
To state them let us fix notations again:
Let $\pi_Y:k\langle\langle X_0, X_1\rangle\rangle\to
k\langle\langle Y_1,Y_2,\dots\rangle\rangle$ be 
the $k$-linear map between non-commutative formal power series rings
that sends all the words ending in $X_0$ to zero and the
word $X_0^{n_m-1}X_1\cdots X_0^{n_1-1}X_1$ ($n_1,\dots,n_m\in\bold N$)
to $(-1)^mY_{n_m}\cdots Y_{n_1}$.
Define the coproduct $\Delta_*$ on $k\langle\langle
Y_1,Y_2,\dots\rangle\rangle$ by
$$\Delta_* Y_n=\sum_{i=0}^n Y_i\otimes Y_{n-i}$$
with $Y_0:=1$.
For $\varphi=\sum_{W:\text{word}} c_W(\varphi) W\in U\frak F_2=
k\langle\langle X_0,X_1\rangle\rangle$
with $c_W(\varphi)\in k$
(a \lq word' means a monic monomial element or $1$ in $U\frak F_2$
%$k\langle\langle X_0,X_1\rangle\rangle$
), 
%define the {\it series shuffle regularization}
put
$$\varphi_*=\exp\left(\sum_{n=1}^{\infty}
\frac{(-1)^n}{n}c_{X_0^{n-1}X_1}(\varphi)Y_1^n\right)\cdot\pi_Y(\varphi).$$
%\end{equation}
The {\it regularised double shuffle relations}
for a group-like series 
$\varphi\in U\frak F_2$
mean %the equations \eqref{shuffle} and
\begin{equation}\label{stuffle}
\Delta_*(\varphi_*)=\varphi_*\widehat\otimes \varphi_*.
\end{equation}

\begin{rem}\label{double shuffle for MZV}
The {\it regularised double shuffle relations} for MZV's
mean the algebraic relations among them
obtained from \eqref{shuffle} and \eqref{stuffle}
for $\varphi=\varPhi_{KZ}$ (cf. \cite{IKZ, R}).
It is also expected that the relations might
produce all algebraic relations among MZV's.
\end{rem}

The following is the simplest example of the relations.
\begin{eg} For $a,b>1$,
\begin{align*}
\zeta(a)\zeta(b)&=
\zeta(a,b)+\zeta(a+b)+\zeta(b,a) \\
&=\sum_{i=0}^{a-1}\binom{b-1+i}{i}\zeta(a-i,b+i)+
\sum_{j=0}^{b-1}\binom{a-1+j}{j}\zeta(b-j,a+j).
\end{align*}
\end{eg}

Our second result here is the implication of
the regularised double shuffle relations from the pentagon equation.

\begin{thm}[\cite{F11}]\label{pentagon-double shuffle}
Let $\varphi=\varphi(X_0,X_1)$ be a %commutator 
group-like element of $U\frak F_2$. % without linear terms.
Suppose that $\varphi$ satisfies the pentagon equation \eqref{pentagon}.
Then  it also satisfies the regularised double shuffle relations \eqref{stuffle}.
\end{thm}

This result attains the final goal of the project posed by Deligne-Terasoma
\cite{T}.
Their idea is to use some convolutions of perverse sheaves,
whereas our proof is to use Chen's bar construction calculus.

\begin{rem}
Our theorem \ref{pentagon-double shuffle}
was extended cyclotomically in \cite{F10b}.
\end{rem}

The following Zagier's relation %among MZV's
which is essential for Brown's proof of theorem \ref{freeness}
might be also one of the most fascinating ones.
The author does not know if it also follows from
our pentagon equation \eqref{pentagon}.
%or the regularised double shuffle relations \eqref{stuffle}.

\begin{thm}[\cite{Z}]\label{Zagier's formula}
$$\zeta(2^{\{a\}},3,2^{\{b\}})=2\sum_{r=1}^{a+b+1}(-1)^r(A^r_{a,b}-B^r_{a,b})
\zeta(2r+1)\zeta(2^{\{a+b+1-r\}})$$
with $A^r_{a,b}=\binom{2r}{2a+2}$
and $B^r_{a,b}=(1-2^{-2r})\binom{2r}{2b+1}$.
\end{thm}

%%%%%%%%%%%%%%%%%%%%%%%%%%%%%%%%%%%%%%%%%%%%%%%%%%%%%%%%%%%%%%%%%%%%%%
\section{Four Groups}
We %review  definitions and 
explain recent developments on the four 
pro-unipotent algebraic groups related to associators;
the motivic Galois group,
the Grothendieck-Teichm\"{u}ller group,
the double shuffle group
and the Kashiwara-Vergne group.

%%%%%%%%%%%%%%%%%%%%%%%%%%%%%%%%%%%%%%%%%%%%%%%%%%%%%%%%%%%%%%%%%%%%%%
\subsection{\bf Motivic Galois group.}
We review on the motivic Galois group, % which is defined to be 
the tannakian dual group
of the category of unramified mixed Tate motives.

We work in the triangulated category $DM(\bold Q)_{\bold Q}$
of {\it mixed motives}
\footnote{
A part of idea of mixed motives is explained  \cite{De} \S 1.
According to {\it Wikipedia},
{\it \lq\lq the (partly conjectural) theory of motives is an attempt to find a universal way to linearize algebraic varieties, i.e. motives are supposed to provide a cohomology theory which embodies all these particular cohomologies."}
}
over $\bold Q$
constructed by Hanamura, Levine and Voevodsky.
%Levine \cite{L2} showed an equivalence of these two categories.
{\it Tate motives} $\bold Q(n)$ ($n\in\bold Z$) are
(Tate) objects of the category.
Let $DMT(\bold Q)_{\bold Q}$ be the triangulated sub-category of 
$DM(\bold Q)_{\bold Q}$ generated by Tate  motives $\bold Q(n)$ ($n\in\bold Z$).
By the work of Levine a neutral tannakian $\bold Q$-category 
$MT(\bold Q)=MT(\bold Q)_{\bold Q}$ of {\it mixed Tate motives over $\bold Q$} 
is extracted by taking a heart with respect to a $t$-structure of
$DMT(\bold Q)_{\bold Q}$.
Deligne and Goncharov \cite{DG} introduced the full subcategory
$MT(\bold Z)=MT(\bold Z)_\bold Q$
of {\it unramified mixed Tate motives} inside. %in $MT(\bold Q)_{\bold Q}$,
All objects  there are mixed Tate motives $M$ (i.e. an object of $MT(\bold Q)$)
such that for each subquotient $E$ of $M$ 
which is an extension of $\bold Q(n)$ by $\bold Q(n+1)$ for $n\in\bold Z$,
the extension class of $E$ in 
$$
Ext^1_{MT(\bold Q)}(\bold Q(n),\bold Q(n+1))=
Ext^1_{MT(\bold Q)}(\bold Q(0),\bold Q(1))=\bold Q^\times\otimes{\bold Q}
$$
is equal to $\bold Z^\times\otimes\bold Q= \{0\}$.

In the category $MT(\bold Z)$ %=MT(\bold Z)_\bold Q$
of unramified mixed Tate motives, the followings hold:
\begin{align}\label{Ext^1}
\dim_{\bf Q} \ & Ext_{MT({\bold Z})}^1({\bold Q}(0),{\bold Q}(m))=
\begin{cases}
1\ (m=3,5,7,\dots),\\
0\ (m:\text{others}),
\end{cases}\\  \label{Ext^2}
\dim_{\bf Q} \ & Ext_{MT({\bold Z})}^2({\bold Q}(0),{\bold Q}(m))=0.\\ \notag
\end{align}

The category $MT(\bold Z)$ 
forms a neutral tannakian $\bold Q$-category
with the fiber functor  $\omega_\mathrm{can}:MT(\bold Z)\to Vect_{\bold Q}$
($Vect_{\bold Q}$: the category of $\bold Q$-vector spaces) 
sending each motive $M$ to
$\oplus_n Hom(\bold Q(n),Gr^W_{-2n}M)$.

\begin{defn}
The {\it motivic Galois group} here %of unramified mixed Tate motives $MT(\bold Z)$ 
is defined to be
the pro-$\bf Q$-algebraic group
$\text{Gal}^\mathcal{M}(\bold Z)
:=\underline{Aut}^\otimes(MT(\bold Z):\omega_\mathrm{can})$.
\end{defn}

By tannakian category theory,  $\omega_\mathrm{can}$ induces an equivalence of categories
\begin{equation}\label{tannakian equivalence}
MT(\bold Z)\simeq Rep \text{Gal}^\mathcal{M}(\bold Z)
\end{equation}
where RHS means the category of finite dimensional $\bold Q$-vector spaces
with $\text{Gal}^\mathcal{M}(\bold Z)$-action.

\begin{rem}
The action of $\text{Gal}^\mathcal{M}(\bold Z)$ 
on $\omega_\mathrm{can}(\bold Q(1))=\bold Q$
defines a surjection $\text{Gal}^\mathcal{M}(\bold Z)\to\bold G_m$ 
and its kernel $\text{Gal}^\mathcal{M}(\bold Z)_1$
is the unipotent radical of $\text{Gal}^\mathcal{M}(\bold Z)$.
There is a canonical splitting $\tau:{\bf G}_m\to \text{Gal}^\mathcal{M}(\bold Z)$
which gives a negative grading on its associated Lie algebra 
$Lie \text{Gal}^\mathcal{M}(\bold Z)_1$.
From \eqref{Ext^1} and \eqref{Ext^2} it follows that
the Lie algebra is the graded {\it free} Lie algebra
generated by one element in each degree $-3,-5,-7,\dots$.
(consult  \cite{De} \S 8 for the full story).
\end{rem}

The {\it motivic fundamental group}
$\pi_1^{\mathcal M}({\bold P}^1\backslash\{0,1,\infty\}
:\overrightarrow{01})$
%with $X=\bold P^1\backslash\{0,1,\infty\}$
constructed in \cite{DG} \S4
is a (pro-)object of $MT(\bold Z)$.
By our tannakian equivalence \eqref{tannakian equivalence},
it gives a  (pro-)object of RHS of \eqref{tannakian equivalence},
which induces a (graded) action
%$\text{Gal}^\mathcal{M}(\bold Z)\to {Aut}{U\frak F_2}$
%where $\underline{F_2}=\omega_\text{can}(\pi_1^{\mathcal M}(X:\overrightarrow{01}))$ is the free pro-unipotent algebraic group of rank $2$.
%Denote its restriction into the unipotent part by
\begin{equation}\label{map}
\Psi:\text{Gal}^\mathcal{M}(\bold Z)_1\to {Aut}\exp{\frak F}_2.
%\to \underline{Aut}\underline{F_2}.
\end{equation}
%Here $\underline{F_2}$ stands for the free pro-unipotent algebraic group
%with two generators $e^{X_0}$ and $e^{X_1}$ and
%$\underline{F_2}(k)$ is equal to the subset of group-like series in $U\frak F_2$.
%By $\underline{Aut}\exp{\frak F_2}$ we mean the pro-algebraic group
%which represents $k\mapsto {Aut}\exp{\frak F_2}(k)$.

\begin{rem}
For each $\sigma\in \text{Gal}^\mathcal{M}(\bold Z)_1(k)$,
its action on $\exp{\frak F_2}$ is described by
$e^{X_0}\mapsto e^{X_0}$ and $e^{X_1}\mapsto 
\varphi_\sigma^{-1} e^{X_1}\varphi_\sigma$
for some $\varphi_\sigma\in \exp{\frak F_2}$.
\end{rem}

The following has been conjectured (Deligne-Ihara conjecture) for a long time
and finally proved by Brown
by using  Zagier's relation (Theorem \ref{Zagier's formula}).

\begin{thm}[\cite{B}]\label{freeness}
The map $\Psi$ is injective.
\end{thm}

It is a unipotent analogue of the so-called Bely\u{\i}'s theorem.
The theorem says that all unramified mixed Tate motives are
associated with MZV's.

%%%%%%%%%%%%%%%%%%%%%%%%%%%%%%%%%%%%%%%%%%%%%%%%%%%%%%%%%%%%%%%%%%%%%%
\subsection{\bf Grothendieck-Teichm\"{u}ller group.}
The Grothendieck-Teichm\"{u}ller group was introduced by Drinfeld \cite{Dr}
in his study of deformations of quasi-triangular quasi-Hopf
quantized universal enveloping algebras.
It was defined to be the set of \lq degenerated' associators.
The construction of the group was also stimulated  by the previous idea of Grothendieck,
{\it un jeu de Teichm\"{u}ller-Lego},
posed in his article
{\it Esquisse d'un programme} \cite{G}.

\begin{defn}[\cite{Dr}]
The %(unipotent part of the graded) 
{\it Grothendieck-Teichm\"{u}ller} %(pro-algebraic) 
{\it group} $GRT_1$
is defined %by ${\underline M}\backslash M$, that is,
to be the pro-algebraic variety 
whose set of $k$-valued points consists of
group-like series $\varphi\in U\frak F_2$ %(i.e. $\varphi\in\underline{F_2}(k)$)
%with $c_{X_0}(\varphi)=c_{X_1}(\varphi)=c_{X_0X_1}(\varphi)=0$
satisfying the defining equations \eqref{pentagon}$\sim$\eqref{hexagon-b}
of associators with $\mu=0$.
\end{defn}

\begin{rem}\label{GRT-group}
(i). By our theorem \ref{pentagon-hexagons}, it is reformulated to be the set of
group-like series satisfying \eqref{pentagon} without quadratic terms.

(ii).
It forms a %pro-unipotent algebraic
group \cite{Dr}
by the multiplication below
\begin{equation}\label{multiplication}
\varphi_2\circ\varphi_1:=\varphi_1(\varphi_2 X_0\varphi_2^{-1},X_1)\cdot\varphi_2
=\varphi_2\cdot\varphi_1(X_0,\varphi^{-1}_2X_1\varphi_2).
\end{equation}

(iii).
By the map sending $X_0\mapsto X_0$ and $X_1\mapsto \varphi^{-1} X_1\varphi$,
the group $GRT_1$ is regarded as a subgroup of ${Aut}\exp{\frak F_2}$.

(iii).
The  cyclotomic analogues of associators and  the Grothendieck-Teichm\"{u}ller group were introduced by Enriquez \cite{E}.
Some elimination results  on their defining equations in special case were obtained 
in \cite{EF}.
\end{rem}

Geometric interpretation (cf. \cite{Dr}) of the equations
\eqref{pentagon}$\sim$ \eqref{hexagon-b} implies the following

\begin{prop}\label{first inclusion}
$\text{Im}\Psi\subset GRT_1$.
\end{prop}

Actually it is expected that they are isomorphic.

\begin{rem}\label{example of GRT}
(i). The Drinfeld associator $\varPhi_{KZ}$ is an associator
(cf. example \ref{Drinfeld associator})
but is not a degenerated associator,
i.e. $\varPhi_{KZ}\not\in GRT_1(\bf C)$.
%but it does not belong to the Grothendieck-Teichm\"{u}ller group $GRT_1$.

(ii). The $p$-adic Drinfeld associator $\varPhi_{KZ}^p$  introduced in \cite{F04}
is not an associator
but a degenerated associator, i.e. $\varPhi_{KZ}^p\in GRT_1({\bf Q}_p)$
(cf. \cite{F07}).
\end{rem}

%%%%%%%%%%%%%%%%%%%%%%%%%%%%%%%%%%%%%%%%%%%%%%%%%%%%%%%%%%%%%%%%%%%%%%
\subsection{\bf Double shuffle group.}
The double shuffle group was introduced by Racinet
as the set of solutions of the regularised double shuffle relations
with \lq degeneration' condition (no quadratic terms condition).

\begin{defn}[\cite{R}]
The {\it double shuffle group} $DMR_0$
is the pro-algebraic variety
whose set of $k$-valued points consists of the group-like series
$\varphi\in U\frak F_2$
%$c_{X_0}(\varphi)=c_{X_1}(\varphi)=c_{X_0 X_1}(\varphi)=0$
satisfying %the generalized double shuffle relation 
the regularised double shuffle relations
\eqref{stuffle}
without linear terms and quadratic terms.
\end{defn}

\begin{rem}
(i).
We note that $DMR$ stands for {\it double m\'{e}lange regularis\'{e}}
(\cite{R}).

(ii).
It was shown  in \cite{R} that it forms a group by 
the operation \eqref{multiplication}.

(iii).
By the same way to remark \ref{GRT-group} (iii),
the group $DMR_0$ is regarded as a subgroup of ${Aut}\exp{\frak F_2}$.
\end{rem}

It is also shown that $\text{Im}\Psi$ is contained in $DMR_0$
(cf.\cite{F07})).
Actually it is expected that they are isomorphic.
Theorem \ref{pentagon-double shuffle} follows
the inclusion between $GRT_1$ and $DMR_0$:

\begin{prop}\label{second inclusion}
$GRT_1\subset DMR_0$.
\end{prop}

It is also expected that they are isomorphic.
%The remarks below are actually  implied by the proposition
%but let us mention  them.

\begin{rem}
(i). The Drinfeld associator $\varPhi_{KZ}$ satisfies
the regularised double shuffle relations
(cf. remark \ref{double shuffle for MZV})
but it is not an element of the double shuffle group,
i.e. $\varPhi_{KZ}\not\in DMR_0(\bf C)$,
%but it does not belong to the Grothendieck-Teichm\"{u}ller group $GRT_1$.
because its quadratic terms are non-zero,
actually $\zeta(2)X_1X_0-\zeta(2)X_0X_1$.

(ii). The $p$-adic Drinfeld associator $\varPhi_{KZ}^p$ % \cite{F04}
satisfies
the regularised double shuffle relations
(cf. \cite{BF, FJ})
and it is an element of the double shuffle group,
i.e. $\varPhi_{KZ}^p\in DMR_0({\bf Q}_p)$,
which also follows from remark \ref{example of GRT}.(ii) and
proposition \ref{second inclusion}.
\end{rem}

%%%%%%%%%%%%%%%%%%%%%%%%%%%%%%%%%%%%%%%%%%%%%%%%%%%%%%%%%%%%%%%%%%%%%%
\subsection{\bf Kashiwara-Vergne group.}
In \cite{KV} Kashiwara and Vergne proposed a conjecture
relating on Campbell-Baker-Hausdorff series which generalises
Duflo's theorem (Duflo isomorphism) to some extent.
%Kashiwara-Vergne conjecture was posed by them in \cite{KV}.
The conjecture was settled generally by Alekseev and Meinrenken \cite{AM}.
The Kashiwara-Vergne group was introduced as a \lq degeneration'  of the set of
solution of the conjecture
by Alekseev and Torossian in \cite{AT},
where they gave another proof of the conjecture
%proved the Kashiwara-Vergne conjecture posed in \cite{KV}
by using associators.

The following is one of %different re
formulations of the conjecture
stated in \cite{AET}.% (as for the original version see \cite{KV})

{\bf Generalized Kashiwara-Vergne problem}:
%(the version in \cite{AET}):
Find a group automorphism $P:\exp\frak F_2\to\exp\frak F_2$
such that %for some $\mu\in k^\times$
$P$ belongs to $TAut\exp\frak F_2$
(that is,
$$P(e^{X_0})=p_1e^{X_0}p_1^{-1} \text{  and  } P(e^{X_1})=p_2 e^{X_1} p_2^{-1}$$ %, \qquad$$
for some $p_1$, $p_2\in\exp\frak F_2$)
and $P$ satisfies
$$
P(e^{X_0}e^{X_1})=e^{(X_0+X_1)}
$$
and the coboundary Jacobian condition %(see \cite{AET})
$$
\delta\circ J(P)=0.
$$

Here
%$F_2(k)$ means the prounipotent free algebraic group with two variables
%$x_0$ and $x_1$. The symbol
$J$ stands for the Jacobian cocycle $J:TAut\exp\frak F_2\to\frak{tr}_2$
and $\delta$ means the differential map $\delta:\frak{tr}_n\to\frak{tr}_{n+1}$
for $n=1,2,\dots$
(for their precise definitions see \cite{AT}).
We note that $P$ is uniquely determined by the pair $(p_1,p_2)$.

The following is essential for the proof of the conjecture.
% in \cite{AT}.

\begin{prop}[\cite{AT, AET}]\label{solution of KV-problem}
Let $(\mu,\varphi)$ be an associator.
Then the pair
$$
(p_1,p_2)=\Bigl(
\varphi(X_0/\mu,X_\infty/\mu), \ e^{X_\infty/2}\varphi(X_1/\mu,X_\infty/\mu)
\Bigr)
$$
with $X_\infty=-X_0-X_1$ gives a solution to the above problem.
\end{prop}

The Kashiwara-Vergne group is %\lq the degeneration' of
defined to be the set of solutions of the problem
with \lq degeneration condition' (\lq the condition of $\mu=0$'):

\begin{defn}[\cite{AT,AET}]
The {\it Kashiwara-Vergne group} $KRV$ is defined to be the set of
$P\in Aut\exp \frak F_2$
%consiting of $P\in Aut\exp \frak F_2$
which satisfies
$P\in TAut\exp\frak F_2$,
%$P(e^{X_0})=p_1e^{X_0}p_1^{-1}$, $P(e^{X_1})=p_2 e^{X_1} p_2^{-1}$,
$$P(e^{(X_0+X_1)})=e^{(X_0+X_1)}$$
and the coboundary Jacobian condition % (see \cite{AET})
$\delta\circ J(P)=0$.
\end{defn}

It forms a subgroup of $Aut\exp \frak F_2$.
We denote by $KRV_{0}$ the subgroup of $KRV$ consisting of 
%automorphism 
$P$ without linear terms in both $p_1$ and $p_2$.
Proposition \ref{solution of KV-problem} yields the inclusion below.

\begin{prop}\label{third inclusion}
$GRT_1\subset KRV_{0}$.
\end{prop}

Actually it is expected that they are isomorphic (cf. \cite{AT}).
Recent result of Schneps in \cite{S} also leads

\begin{prop}\label{fourth inclusion}
$DMR_0\subset KRV_{0}$.
\end{prop}

\subsection{Comparison}
By combining theorem \ref{freeness} and 
proposition \ref{first inclusion}, \ref{second inclusion},
\ref{third inclusion} and \ref{fourth inclusion},
we obtain
\begin{prop}
$
\text{Gal}^\mathcal{M}(\bold Z)_1
\subseteq GRT_1
\subseteq DMR_0
\subseteq KRV_0
$.
\end{prop}

%Now it would be enough interesting to ask
Here we come to an interesting question on our four groups.
\begin{q}
Are they all equal? Namely,
$$
\text{Gal}^\mathcal{M}(\bold Z)_1=GRT_1=DMR_0=KRV_0 \ \ ?
$$
\end{q}

Though it might be not so good mathematically to believe such equalities without
a  strong conceptual support,
%believe the equality.
%However 
the author thinks that 
it might be good at least spiritually to imagine/expect
a hidden theory to relate them behind.

%%%%%%%%%%%%%%%%%%%%%%%%%%%%%%%%%%%%%%%%%%%%%%%%%%%%%%%%%%%%%%%%%%%%%%

\end{document}